\documentclass[12pt]{amsart}
\usepackage{amssymb}
\usepackage{verbatim}

\theoremstyle{plain}
\newtheorem{theorem}{Theorem}
\newtheorem{corollary}{Corollary}
\newtheorem{proposition}{Proposition}
\newtheorem{lemma}{Lemma}

\theoremstyle{definition}
\newtheorem*{definition}{Definition}
\newtheorem{example}{Example}

\theoremstyle{remark}
\newtheorem*{remark}{Remark}

\def\bbC{\mathbb C}

\def\bbN{\mathbb N}

\def\bbZ{\mathbb Z}

     \newcommand{\sA}{\mathcal A}
     \newcommand{\sB}{\mathcal B}
     \newcommand{\sC}{\mathcal C}
     \newcommand{\sD}{\mathcal D}
     \newcommand{\sE}{\mathcal E}
     \newcommand{\sF}{\mathcal F}
     \newcommand{\sG}{\mathcal G}
     \newcommand{\sH}{\mathcal H}
     
     \newcommand{\sK}{\mathcal K}
     
     \newcommand{\sL}{\mathcal L}
     
     \newcommand{\sN}{\mathcal N}

     \newcommand{\sS}{\mathcal S}
     \newcommand{\sT}{\mathcal T}

     \newcommand{\sX}{\mathcal X}

\newcommand{\cstar}{\ensuremath{\text{C}^{*}}-}

\DeclareMathOperator{\dist}{dist}
\DeclareMathOperator{\supp}{supp}

\DeclareMathOperator{\spec}{spec}
\DeclareMathOperator{\Alg}{Alg}
\DeclareMathOperator{\Lat}{Lat}
\DeclareMathOperator{\Lie}{Lie}
\DeclareMathOperator{\tr}{tr}
\newcommand{\ds}{\displaystyle}

\begin{document}

% topmatter

\title[Lie Ideals]{Lie Ideals in Operator Algebras}
\author{Alan Hopenwasser}
\address{Department of Mathematics\\
	University of Alabama\\
	Tuscaloosa, AL 35487-0350}
\email{ahopenwa@euler.math.ua.edu}

 \author{Vern Paulsen}
 \address{Department of Mathematics\\
	University of Houston\\
	Houston, TX 77204-3476}
 \email{vern@math.uh.edu}

\thanks{This research is supported in part by a grant from the 
NSF}
 \thanks{The authors  thank Allan Donsig, Robert 
L. Moore,
   and David Pitts for helpful comments on the material of this 
paper.  The authors also thank Ken Davidson for organizing a 
workshop
on non-self-adjoint algebras at the Fields Institute in July 
2002.
Substantial improvements to the original version of this paper 
were
made at the workshop.}
 \date{\today}

% abstract
 \begin{abstract}
Let $\mathcal A$ be a Banach algebra for which the group of 
invertible
elements is connected.  A subspace 
$\mathcal L \subseteq \mathcal A$ is a
Lie ideal in $\mathcal A$ if, and only if, it is invariant under 
inner
automorphisms.  This applies, in particular, to any canonical
subalgebra of an AF \ensuremath{\text{C}^{*}}-algebra.  The same
theorem is also proven for strongly closed subspaces of a totally
atomic nest algebra whose atoms are ordered as a subset of the
integers and for CSL subalgebras of such nest algebras.

We also give a detailed description of the structure of a Lie 
ideal in
any canonical triangular subalgebra of an AF 
\ensuremath{\text{C}^{*}}-algebra.
 \end{abstract}
\maketitle

\section{Introduction} \label{intro}

In view of the close relationship between derivations and
automorphisms, it is not surprising that in many settings a 
subspace
of an algebra is a Lie ideal if, and only if, it is invariant 
under
similarity transformations.  We prove  this equivalence for 
closed
subspaces of any Banach algebra for which the group of invertible
elements is connected.  This includes all canonical subalgebras 
of an
AF \cstar algebra.  The proof of this result is short and direct.
Initially, we proved the equivalence in the context of triangular
subalgebras of AF \cstar algebras via a detailed analysis of the
structure of Lie ideals in triangular subalgebras.  While the
structure theorem is no longer needed to prove that Lie ideals 
are
similarity invariant, it remains of independent interest and is
described in section \ref{AF} of this paper.

In the process of investigating triangular subalgebras of 
AF \cstar algebras, it is appropriate to look at triangular
subalgebras of finite dimensional \cstar algebras.  In fact, with
only a moderate
 additional effort, we can obtain a description of Lie ideals
in an arbitrary digraph
algebra.  In all likelihood
these finite dimensional results are  not new, but
the authors know of no suitable reference (except in more
 specialized contexts).  These results appear in section
\ref{digraph}.

The impetus for this note comes from a similar result by Marcoux 
and
Sourour \cite{MR1883197} in a much more limited context: direct 
limits of
full upper triangular matrix algebras ($T_n$'s); i.e., 
subalgebras of
UHF \cstar algebras which are strongly maximal triangular in 
factors.
  The direct limit algebra
context constitutes only a small portion of  \cite{MR1883197}; 
most of
that paper is devoted to weakly closed Lie ideals in nest 
algebras and
to Lie ideals in algebras of infinite multiplicity.
In section \ref{Lie_sim}, where we present the main theorem, we 
also
prove that strongly closed Lie ideals are similarity invariant 
in the
context of totally atomic nest algebras whose atoms are ordered 
as a
subset of the integers.  Since weak and strong closure are 
identical
for subspaces, this result is contained in \cite{MR1883197}.  Our
proof is much shorter than the one in \cite{MR1883197}, at the 
price
of omitting a considerable amount of information about the 
structure
of Lie ideals in nest algebras.  On the other hand, our method 
also
works for CSL-subalgebras of these ``integer-ordered'' nest 
algebras,
so the domain of validity of the  equivalence is extended.

If $\sA$ is an algebra, a subspace $\sL$ is a \emph{Lie ideal} if
$[x,a]=xa-ax \in \sL$ whenever $x \in \sL$ and $a \in \sA$.  The
subspace $\sL$ is said to be \emph{similarity invariant} if
$t^{-1}xt \in \sL$ whenever $x \in \sL$ and $t$ is an invertible
element of $\sA$.  David Pitts has pointed out to the authors
an attractive reformulation of the equivalence of these two 
concepts
(when valid): the family of inner derivations of $\sA$ and the
family of inner automorphisms of $\sA$
have the same invariant subspaces.

In order to avoid any ambiguity in the sequel, we shall use the 
term
``associative ideal,'' rather than the usual term ``ideal,'' for
an ordinary (2-sided) ideal.  Thus, all associative ideals are 
also
Lie ideals, but not conversely.

For closed subspaces of a Banach algebra, similarity invariance 
for
a subspace implies that the subspace is a Lie ideal.  This is an
unpublished result of Topping; a brief proof is
contained in Theorem \ref{Th_sim} of this paper.
 The description of Lie ideals goes back a 
long way; in
a purely algebraic context Herstein \cite{MR16:789e}
studied Lie ideals (and
their relationship with associative ideals) in 1955.  An 
extensive
treatment of the algebraic theory appears in his book
\cite{MR42:6018}.
 Lie ideals in
the algebra of all linear transformations on an infinite 
dimensional
vector space were studied by Stewart in \cite{MR47:6800}. Murphy
\cite{MR85b:16030} investigated
 Lie ideals and their relationship with
associative ideals in algebras with a set of $2 \times 2$ matrix
units. Fong, Meiers and Sourour \cite{MR84g:47039} and Fong and 
Murphy \cite{MR89k:47067} have written about these ideas in the
$\sB(\sH)$ context.  Marcoux \cite{MR96f:46123} identified all 
Lie
ideals in a UHF \cstar algebra and proved that they are 
similarity
invariant (as well as invariant under unitary conjugation).  He 
also
described the Lie ideals in algebras of the form
$\sA \otimes C(X)$, where $\sA$ is either a full matrix algebra 
or a
UHF \cstar algebra.  Further relevant information in the \cstar
algebra setting can be found in Pedersen \cite{MR82b:46078} and 
in
Marcoux and Murphy \cite{MR99b:46085}.  In \cstar algebra 
contexts,
invariance under unitary conjugation is generally equivalent to
invariance under inner derivations.  Moving to  the 
non-self-adjoint
operator algebra literature, see Hudson, Marcoux and Sourour
\cite{MR98k:47090} for a description of the form of Lie ideals 
in nest
algebras and in direct limit algebras which are strongly maximal
triangular in factors.  And, as mentioned above, 
\cite{MR1883197} shows
that a weakly closed subspace in a nest algebra is a Lie ideal 
if, and
only if, it is similarity invariant.  The assumption of weak 
closure
can be dropped if the nest has no finite dimensional atoms.

\section{Lie Spaces and Similarity} \label{Lie_sim}
We begin with a result that refines the relationship between Lie 
ideals and 
similarity invariant subspaces given by Topping.

If $\sA$ is a unital Banach algebra and $\sX$ is a Banach 
space, then we shall call $\sX$ a \emph{bounded, Banach 
$\sA$-bimodule} provided that $\sX$ is an $\sA$-bimodule such 
that the identity element of $\sA$ acts as the identity on 
$\sX$ and provided that the module action is bounded; that is,
there exists a constant $K$ such that $\|axb\| \le 
K\|a\|\|x\|\|b\|$ for all $a,b$ in $\sA$ and $x$ in $\sX$.
A linear subspace(not necessarily a submodule!) $\sL$ of $\sX$ 
is called a \emph{Lie subspace} over $\sA$ provided that 
$ax -xa \in \sL$ for every $x \in \sL$ 
and every $a \in \sA$.
Thus, a Lie ideal is just a Lie subspace of $\sA$.
We call a subspace $\sL$ of $\sX$  \emph{similarity invariant} 
provided that $a^{-1}xa \in \sL$ for every $x \in \sL$ and 
every invertible element $a \in \sA$.
 
\begin{theorem} \label{Th_sim}
Let $\sA$ be a unital Banach algebra, let $\sX$ be a bounded, 
Banach $\sA$-bimodule, let $\sG$ denote the 
connected component
of the identity in the group of invertible
elements of $\sA$ and let $\sL$ be a closed subspace of $\sX$.  
Then
$\sL$ is a Lie subspace if, and only if, $b^{-1}\sL b \subseteq 
\sL$ for every $b \in \sG.$ 
\end{theorem}

\begin{proof}
Assume that $\sL$ is a closed Lie subspace and that $b$ is in 
$\sG.$
  Since $b$ is in the connected component of the identity, $b$
is a finite product of exponentials.
Therefore, to prove that $b^{-1}\sL b \subseteq \sL,$ it 
suffices to prove that
$e^{-a}\sL e^a \subseteq \sL$, for any $a \in \sA$.

Fix $x\in \sL$ and set
$x(t) = e^{-ta}xe^{ta}$.  This is an analytic function.  An easy
induction argument shows that, for all $n \geq 0$, 
the derivatives satisfy the relation
$x^{(n+1)}(t) = x^{(n)}(t)a - a x^{(n)}(t)$.  Since
$x(0)=x \in \sL$, it follows that 
$x^{(n)}(0) \in \sL$  for all $n \geq 0$.  Therefore, all the 
terms in
the power series for $x(t)$ lie in $\sL$. Since $\sL$ is closed, 
it
follows that $x(t) \in \sL$ for all $t$.  In particular,
$e^{-a}xe^a \in \sL$.

Conversely, assume that $b^{-1} \sL b \subseteq \sL$ for every 
$b \in \sG.$
Given any $a \in \sA,$ form $x(t)$ as above. By assumption, 
$x(t) \in \sL$ for all $t$ and hence
 the derivative $x'(t) \in \sL$ for all $t$. Evaluating at 
$t=0$, we find that $xa - ax \in \sL$ and our proof
 is complete.
 \end{proof}

As mentioned in the introduction, Topping has proven that any 
closed
subspace of $\sA$ which is similarity invariant is a Lie ideal. 
His proof 
is essentially
reproduced in the proof of the converse in Theorem \ref{Th_sim}.

Theorem \ref{Th_sim} shows that in order 
to determine whether or not a 
closed Lie subspace $\sL$ of a bounded, Banach $\sA$-bimodule
is similarity invariant, it is sufficient to check whether 
or not $b^{-1} \sL b \subseteq
\sL$ for any collection of elements $b$ that contains at least 
one representative from each coset
in $\sA^{-1}/\sG.$

\begin{corollary}
Let $\sB$ be an AF \cstar algebra with canonical masa $\sD$ and 
let
$\sA$ be a canonical subalgebra of $\sB$, i.e., a subalgebra 
such that
$\sD \subseteq \sA \subseteq \sB$.  A closed subspace of $\sB$ 
is a
Lie subspace over $\sA$ if, and only if, it is invariant under
similarities. 
\end{corollary}

\begin{proof}
Clearly, $\sB$ is a bounded, Banach $\sA$-bimodule.
The invertibles in a canonical subalgebra are connected.  (Any
invertible $t$ can be closely approximated by -- and hence path
connected to -- an invertible in a finite dimensional 
approximant of
$\sA$.  Each invertible in a (finite dimensional) digraph 
algebra is
path connected to the identity element.)
\end{proof}

We now turn attention to some atomic nest algebras.
It is not known whether the invertibles in a nest algebra are
connected, not even when the nest is atomic.  Therefore 
Theorem \ref{Th_sim} does not apply.  However, in the ``integer
ordered'' cases we can still obtain the similarity invariance of 
Lie
ideals without using the structure of Lie ideals. 
Thus we provide, albeit only in a special case, a shortcut to the
argument in  \cite{MR1883197}.
 This method also
works for certain CSL subalgebras of such a nest algebra.

Let $N$ be a subset of $\bbZ$ and, for each $n \in N$, let
$\sH_n$ be a Hilbert space.  When $N$ is a finite set, the 
following
discussion is valid with some minor modification.  It is, 
however,
easy to provide an even simpler proof of Theorem \ref{nest} for 
finite
nests.  Accordingly, we assume that $N$ is an infinite set.  
Without
any loss of generality, we may assume that $N$ is one of
$\bbZ$, $\bbN$ or $-\bbN$.  Let
$\sH = \sum_{n \in N}^{\oplus} \sH_n$.  For each $n$, let $E_n$ 
denote
the orthogonal projection of $\sH$ onto $\sH_n$.  If 
$P_n = \bigvee_{k \leq n}E_n$, then
$\sN = \{P_n \mid n \in N\} \cup \{0,I\}$ is a totally atomic 
nest in
$\sH$ whose atoms, $\{E_n\}_{n\in N}$, are order isomorphic to 
$N$.  ($E_n \ll E_m$ if, and only if, 
$E_n \sH E_m \subseteq \Alg \sN$.)

Let $\sA$ be a reflexive subalgebra of $\Alg \sN$ such that
$\Lat \sA$ is a totally atomic lattice whose atoms are exactly 
the
atoms of $\sN$.  The elements of $\Alg \sN$ consist of all upper
triangular matrices with respect to the decomposition
$\sH = \sum_{n \in N}^{\oplus} \sH_n$.  The elements of $\sA$ 
consist
of those matrices in $\Alg \sN$ whose entries are 0 in certain
specified locations.

Let $A = (A_{i,j})$ be an operator in $\sA$. 
 Each entry $A_{i,j}$ is an
operator in $\sB(\sH_j,\sH_i)$ and $A_{ij} = 0$ when
$i > j$ and when $(i,j)$ is  one of the specified locations 
mentioned
above.  
For each $n \in \bbN \cup \{0\}$, define
$D_n = \sum_{k \in N} E_k A E_{k+n}$.  The matrix for $D_n$ is
$(C_{i,j})$, where $C_{i,i+n} = A_{i,i+n}$ for all $i$ and 
$C_{i,j} = 0$ for all other values of $i$ and $j$.  Now define
\[
A(z) = \sum_{n=0}^{\infty} D_n z^n = \left( 
A_{i,j}z^{j-i}\right).
\]
Note that 
$\|D_n\| \leq \sup_{k \in N} \|E_k A E_{k+n}\| \leq \|A\|$, for 
each
$n \geq 0$. Consequently, the series 
$\sum_{n=0}^{\infty} D_n z^n$ converges uniformly on any disk
$|z| < r <1$ and $A(z)$ is analytic on the open disk $|z|<1$.

If $|z|=1$, then $A(z)$ is unitarily equivalent to $A$.  Indeed, 
write
$z = e^{i\theta}$ and let $U(\theta)$ be the diagonal unitary 
matrix
whose $n^{\text{th}}$-diagonal entry is $e^{in\theta}E_n$.
Then $U(\theta)\in \sA$ and
$U(\theta)^*\sA U(\theta) = A(e^{i\theta})$.  Thus
$\|A(e^{i\theta})\| = \|\sA\|$ for all $\theta$; by the maximum
modulus principle, $\|\sA(z)\| \leq \|\sA\|$ for all
$|z| \leq 1$.  Although the series 
$\sum_{n=0}^{\infty} D_n z^n$ need not converge uniformly on the 
whole
unit disk, it does converge strongly.  This follows from the 
fact that
for any fixed vector $h \in \sH$, $\|D_n h\| \to 0$.  The 
function
$A(z)$ is continuous with respect to the strong operator 
topology on
the closed unit disk.  For any vectors $h_1$ and $h_2$ in $\sH$, 
the
function $z \to \langle A(z)h_1, h_2 \rangle$ is a complex valued
 analytic function in the open unit disk with continuous boundary
 values. 

Finally, observe that if $A \in \sA$ is invertible with inverse 
$B$ in
$\sA$, then $A(z)B(z)=I$ for all $|z| \leq 1$.  Indeed, if
$z = e^{i\theta}$ then
\begin{align*}
A(z)B(z) &= U(\theta)^* A U(\theta)  U(\theta)^* B U(\theta) \\
 &= U(\theta)^* AB U(\theta) =I. 
\end{align*}
Since this identity holds on the boundary of the unit disk and
$A(z)B(z)$ is analytic, it holds throughout the unit disk.

\begin{theorem} \label{nest}
Let $\sA \subseteq B(\sH)$ be a CSL subalgebra of a nest 
algebra $\Alg \sN$ whose 
atoms
have order type isomorphic to a subset of the integers.  Assume 
that
$\Lat \sA$ is totally atomic and that the atoms for $\Lat \sA$ 
are
precisely the atoms for $\sN$.  Let $\sL \subseteq B(\sH)$ be a 
strongly closed 
Lie
subspace over $\sA$.  Then $\sL$ is invariant under 
similarities from 
$\sA$.
\end{theorem}

\begin{proof}
Let $X \in \sL$ and let $A$ be an invertible element of $\sA$ 
with
inverse $B$.  For $|z|<1$, it is easy to see that $A(z)$ is in 
the
connected component of the identity in the group of invertibles 
for
$\sA$.  Since $B(z)$ is the inverse of $A(z)$, 
Theorem \ref{Th_sim} implies that 
$A(z)X B(z) \in \sL$ for all $|z| < 1$.  But $A(z) \to A$ and
$B(z) \to B$ strongly as $z \to 1$ and both $A(z)$ and $B(z)$ are
uniformly bounded on the unit disk, so
$A(z) X B(z) \to AXB$ strongly.  Since $\sL$ is strongly closed,
$AXB = AXA^{-1} \in \sL$.
\end{proof}

\section{Digraph Algebras} \label{digraph}

In this section, we shall describe all the Lie ideals in a 
family of
operator algebras known variously as ``digraph algebras,'' 
``incidence
algebras'' and ``finite dimensional CSL-algebras.''  In 
addition, we
shall give an alternate proof
 that every Lie ideal is similarity invariant. (Since the 
invertibles
 are connected in a digraph algebra, this result is a special 
case of
 Theorem \ref{Th_sim}.)

Fix a finite dimensional Hilbert space $\sH$.  A digraph algebra 
is a
subalgebra $\sA$
of $\sB(\sH)$ which contains a maximal abelian self-adjoint
subalgebra $\sD$ of $\sB(\sH)$.  Since $\sD$ is maximal abelian, 
the
invariant projections for $\sA$, $\Lat \sA$, are elements of 
$\sD$ and
so are mutually commuting.  Thus $\sA$ is a CSL-algebra.  
Obviously,
$\sA$ is finite dimensional; on the other hand,
every finite dimensional CSL-algebra acts
on a finite dimensional Hilbert space and contains a masa.
For another description of $\sA$, let $n$ be the dimension of 
$\sH$.
Then $\sA$ is isomorphic to a subalgebra of $M_n$ which contains 
all
the diagonal matrices, $\sD_n$.  An $n \times n$ pattern
matrix, whose entries consist of $0$'s and $*$'s, is associated 
with $\sA$.  After identifying
$\sA$ with the matrix algebra to which it is isomorphic,  $\sA$
consists of all those matrices with arbitrary entries where 
there are
$*$'s in the pattern matrix
 and $0$'s in the remaining locations.  Not
every pattern gives rise to an algebra, but those that do yield 
all
the digraph algebras.  This description is the one which gives 
rise to
the term ``incidence algebra.''  The term digraph algebra refers 
to
the fact that associated with $\sA$ there is a directed graph on 
the
set of vertices $\{1,2,\dots,n\}$.  This graph contains all the 
self
loops.  Then $\sA$ contains the matrix unit $e_{ij}$
if, and only if, there is a (directed) edge from $j$ to $i$ in 
the
digraph. The matrix units in
 $\sA$ generate $\sA$ as an algebra.

Of these three descriptions, we shall primarily use the incidence
algebra pattern.  Furthermore, for a suitable choice of matrix 
units,
we may assume that $\sA$ has a block upper triangular format.  
One way
to see this is to look at the set $\{f_1, \dots f_p\}$ of minimal
central projections in $\sA \cap \sA^*$.
  It is then easy to show that for each
$i$, $f_i\sA f_i$ is isomorphic to a full matrix algebra and 
that if 
$i \ne j$, then at least one of $f_i \sA f_j$ and $f_j \sA f_i$ 
is 
$\{0\}$.  Furthermore, if $f_i \sA f_j$ contains non-zero 
elements,
then it contains all elements of $f_i \sB(\sH) f_j$.  After a 
possible
reindexing, we may assume that $i >j$ implies $f_i\sA f_j 
=\{0\}$.  A
selection of matrix units for $\sB(\sH)$ compatible with the 
minimal
central projections puts
$\sA$ into block upper triangular form.  An alternate way to 
achieve
the same form is to select a maximal nest from within $\Lat \sA$ 
and
then choose matrix units compatible with
 the nest.  It is then routine to
show that $\sA$ has a block upper triangular form in which each
non-zero block is full.

We shall refer to $\ds \sE = \sum_{i=1}^p f_i\sA f_i$ as the 
diagonal
part of $\sA$ and $\ds \sS = \sum_{i<j} f_i\sA f_j$ as the
off-diagonal part of $\sA$. The diagonal part of $\sA$ contains, 
but
is in general larger than, the masa $\sD$.  Since
$\sE = \sA \cap \sA^*$, the diagonal part is intrinsically
determined.  In contrast, $\sD$ is determined only up to an inner
automorphism from $\sE$.

 Any associative ideal in $\sA$ is, of
course, a Lie ideal; we shall be concerned with associative 
ideals
which are subsets of $\sS$.  One reason for this is that, as we
shall see later, if $\sL$ is a Lie ideal, then
$\sL \cap \sS$ is an associative ideal.

\begin{definition}
An \emph{off-diagonal} associative ideal is an associative ideal 
$\sK$
which is a subset of $\sS$.  This is equivalent to requiring that
$\sK \cap \sD = \{0\}$.
\end{definition}

Fix, for the moment, an off-diagonal associative ideal $\sK$.  
Then
$\sK$ is the smallest Lie ideal $\sL$
with the property that
$\sL \cap \sS = \sK$.  It is a simple matter to check that if
$f_i \sK f_j \ne \{0\}$, then 
$f_i \sK f_j = f_i \sA f_j \;(= f_i \sB(\sH) f_j)$.  In 
addition, for 
$t \leq i$ and $s \geq j$, 
$f_t \sK f_s = f_t \sA f_s$.  Thus, $\sK$ consists only of full
blocks and, when the pattern for $\sK$ contains a $*$,  there is
also a $*$ in all locations above and to the right (based on the
pattern for $\sA$).  Note: we use the non-strict interpretations 
for
``above'' and ``to the right.''  ``Above'' permits entries in 
the same
row and ``to the right'' permits entries in the same column.
What we have just described for off-diagonal ideals is, of 
course, 
true for all associative ideals.

There are Lie ideals larger than $\sK$ whose off-diagonal part
(intersection with $\sS$) is $\sK$.  Each of these can be 
obtained by
adding an appropriate subspace of $\sE$ to $\sK$.  Accordingly, 
we
make the following definition:

\begin{definition}
Let $\sK$ be an off-diagonal associative ideal in $\sA$.  A 
\emph{Lie addend} for $\sK$ is a subspace $\sG$ of $\sE$ with the
property that $\sG + \sK$ is a Lie ideal.
\end{definition}

\begin{example} \label{ex:maxadd}
We describe an example of a Lie addend $\sF$
for $\sK$.  Later on, we shall
see that this is the largest Lie addend for $\sK$; equivalently, 
$\sF + \sK$ is the largest Lie ideal which satisfies the 
property that
$\sL \cap \sS = \sK$.
Let 
$\tilde S = \{(i,j) \mid f_i \sA f_j \ne \{0\} \text{ and } 
i<j\}$ 
and
 $\tilde K = \{(i,j) \mid f_i \sK f_j \ne \{0\}\}$. 
Each $f_i \sE f_i$ is isomorphic to a full matrix algebra. 
Suppose that $(i,j) \in \tilde S \setminus \tilde K$.  Then for 
each
 $x \in \sF$,  both $f_ixf_i$ and $f_jxf_j$ must be scalar 
matrices
 with equal scalars.  On the other hand, if $i$ is such that
no ordered pair $(i,j)$ or $(j,i)$ lies in
$\tilde S \setminus \tilde K$, then $f_ixf_i$ is arbitrary.  In 
this
case, $f_i \sF f_i = f_i \sB(\sH) f_i$ is a subset of $\sF$.
Thus, we see that each $f_i \sF f_i$ is either the scalars or 
the full
matrix algebra $f_i \sB(\sH) f_i$.  There are constraints 
relating
some of the scalar blocks, as indicated above; there are no
constraints involving the full matrix algebra blocks.

To show that $\sF + \sK$ is a Lie ideal, it is
sufficient to show that for each $x \in \sF + \sK$ and each 
matrix
unit $e_{ts} \in \sA$, $[x,e_{ts}] \in \sF + \sK$.  Since $\sK$ 
is a
Lie ideal in $\sA$, we may  restrict attention to the case in
which $x \in \sF$.  The matrix unit $e_{ts}$ is either in $\sS$ 
or in
$\sE$.  First consider the case in which $e_{ts}$ is in $\sS$.  
Then
there is a pair $(i,j) \in \tilde{S}$
 such that $e_{ts} = f_i e_{ts} f_j$. (Note:
 $i<j$).  Now there are
two subcases to consider.  One is when $(i,j) \in \tilde K$.  
Write 
$x = \sum_k f_k x f_k$.  If $k \ne i,j$ then 
$[f_k x f_k, e_{ts}]=0$, an element of $\sF+\sK$.  If $k =i$ or
$k=j$, then  $[f_k x f_k, e_{ts}]\in f_i \sA f_j = f_i \sK f_j$ 
and so is an
element of $\sK$.  It now follows that $[x,e_{ts}] \in \sK$.
  The remaining subcase occurs when
$(i,j) \in \tilde S \setminus \tilde K$.  Then,
$[x,e_{ts}]=f_i x f_i e_{ts}-e_{ts}f_j x f_j$.  But now there is 
a
scalar $\lambda$ such that $f_ixf_i = \lambda f_i$ and
$f_jxf_j = \lambda f_j$; consequently, $[x,e_{ts}]=0$.
The second case to consider is when $e_{ts} \in \sE$.  Then for 
some
$i$, $e_{ts}=f_ie_{ts}f_i$.  It follows that $[x,e_{ts}]$ is an
element of $f_i \sA f_i$ and 
$[x,e_{ts}]=[f_ixf_i,e_{ts}]$.  Either $f_i \sA f_i$ is a 
subset of $\sF$ (when there is no $j$ such that either $(i,j)$ or
$(j,i)$ lies in $\tilde S \setminus \tilde K$)
 or $f_i x f_i$ is scalar and $[x,e_{ts}]=0$; either way,
$[x,e_{ts}] \in \sF + \sK$.
\end{example}

Before describing the structure of Lie addends, we show that any 
Lie
ideal has the form $\sL = \sG + \sK$, where $\sK$ is an 
off-diagonal
associative ideal and $\sG$ is a Lie addend for $\sK$.  

\begin{proposition} \label{dig_Lieid_struct}
Let $\sL$ be a Lie ideal in $\sA$.  Let $\sK = \sL \cap \sS$ and
$\sG = \sL \cap \sE$.  Then
$\sL = \sG + \sK$ and $\sK$ is an associative ideal in
$\sA$.
\end{proposition}

\begin{proof}
As before,
 $f_1, \dots, f_p$ are the minimal central projections of $\sE$ 
in
an order which renders $\sA$ block upper triangular after the
selection of a system of matrix units compatible with the $f_i$.
Define $\pi \colon \sA \to \sE$ by 
$\pi(x) = \sum_i f_i x f_i$. Note that
 $\pi$ is a conditional expectation onto $\sE$. 

Suppose that $x \in \sL$.
We claim that $x - \pi(x)$ is an element 
of  $\sL$.  Since
$\ds x - \pi(x) = \sum_{i<j} f_i x f_j$, it will suffice to show 
that
each $f_i x f_j \in \sL$.  But this follows from the fact that
$f_i x f_j = [f_i, [x,f_j]]$ (since $f_j x f_i = 0$ when $i<j$).
Since both $x$ and $x - \pi(x)$ are in $\sL$, so is $\pi(x)$.
We now have
\[
\sK = \{x \in \sL \mid \pi(x) = 0 \}
=\{x-\pi(x) \mid x \in \sL \}
\]
and
\[
\sG = \pi(\sL) \subseteq \sL.
\]
For any $x \in \sL$, $x = \pi(x) + (x - \pi(x))$, so
$\sL = \sG + \sK$.
Since $\sS$  and $\sL$ are  Lie ideals in $\sA$, 
$\sK$ is also a Lie
ideal. It remains to show that $\sK$ is an associative ideal.

Let $x \in \sK$.  We have just seen that
 $f_i x f_j \in \sL$ whenever $i < j$.
 Since $f_i x f_j$ is clearly in $\sS$, it
is an element of  $\sK$.
 If $f_i \sK f_j \ne \{0\}$, then we
can find $x \in \sK$ such that $f_i x f_j \ne 0$.  If $y$ is any 
rank
one element of $\sA$ such that $f_i y f_i = y$, then
$y f_i x f_j = [y, f_i x f_j]$ is also in $\sK$.  Similarly, if 
$z$ is
any rank one element of $\sA$ such that
$f_j z f_j = z$, then
$y f_i x f_j z = [yf_ixf_j, z]$ is in $\sK$.  But suitable 
choices of
$y$ and $z$ produce any rank one element of $f_i \sB(\sH) f_j$.  
Thus,
$f_i \sK f_j \ne \{0\}$ implies that
$f_i \sB(\sH) f_j \subseteq \sK$; in other words, $\sK$ consists 
of
certain of the strictly upper triangular blocks from $\sA$.

To show that $\sK$ is an associative
ideal in $\sA$, we need to show that when a 
block appears in $\sK$, so does each block (from the pattern for
$\sA$) which lies to the right and above.  So, assume 
  that $f_i \sB(\sH) f_j \subseteq \sK$ 
and (with $j<k$) that $f_j \sA f_k \ne \{0\}$.
Let $y \in f_j \sA f_k$ be such that $y \ne 0$ and let 
$x \in f_i \sB (\sH) f_j$ be such that $xy \ne 0$.  Now
$yx = yf_k f_i x = 0$ (since $i < j < k$), so
$xy = [x,y]$ is a non-zero element of  $\sK$.  
But this implies that $f_i \sB(\sH) f_k \subseteq \sK$.  
In a similar way, if $k < i < j$ and $f_k \sA f_i \ne \{0\}$, 
then
$f_i \sB(\sH) f_j \subseteq \sK$.  Thus $\sK$ is an 
associative ideal in $\sA$.
\end{proof}

To complete the description of Lie ideals in $\sA$, it remains to
describe the structure of an arbitrary Lie addend $\sG$ for an
off-diagonal associative ideal $\sK$.  We continue to identify
$f_i \sB(\sH)f_i$ with a full matrix algebra acting on
$f_i \sH$.

\begin{proposition} \label{Lieadd_struct}
Let $\sK$ be a an off-diagonal associative ideal and let $\sG$ 
be a
Lie addend for $\sK$. 
Let 
$\tilde S = \{(i,j) \mid f_i \sA f_j \ne \{0\} \text{ and } 
i<j\}$ 
and
 $\tilde K = \{(i,j) \mid f_i \sK f_j \ne \{0\}\}$. 
If $(i,j) \in \tilde S \setminus \tilde K$, then for each 
$x \in \sG$, $f_i x f_i$ and $f_j x f_j$ are scalar matrices with
equal scalars.  For each $i$, $f_i \sG f_i$ is one of the 
following
four subspaces of 
$f_i \sB(\sH) f_i$: $\{0\}$, the scalar matrices,
the set of all matrices with trace zero, or the full matrix 
algebra.
Each element of $f_i \sG f_i$ with trace zero is an element of 
$\sG$.

In particular, $\sF$, the Lie addend in example \ref{ex:maxadd} 
is a
maximal Lie addend.  Any subspace $\sG$ of $\sF$ which satisfies 
these
properties is a Lie addend.
\end{proposition}

\begin{remark}
It is not true that $f_i \sG f_i \subseteq \sG$, when $\sG$ is a 
Lie
addend. (This would violate the scalar constraints when 
$(i,j) \in \tilde S \setminus \tilde K$.)  The condition that 
trace
zero elements in $f_i \sG f_i$ are in $\sG$ essentially says that
$\sG$ satisfies no non-scalar constraints.  (Additional scalar
constraints beyond the ones required for membership in $\sF$ may 
be
satisfied.) 
\end{remark}

\begin{proof}
If $\sG$ is a subspace of $\sF$ satisfying the conditions in 
Proposition \ref{Lieadd_struct},
 then a slight variation of the argument in example
\ref{ex:maxadd} showing that $\sF+\sK$ is a Lie ideal shows that
$\sG+\sK$ is a Lie ideal.  The assumption that trace zero 
elements of
$f_i \sG f_i$ are in $\sG$ is used when $x \in \sG$ and 
$f_i x f_i$ is non-scalar.  In this situation, there are matrix 
units
supported in the block associated with $f_i$ such that
$[x,e] \neq 0$.  But when $[x,e] \neq 0$ and $e = f_i e f_i$, 
$[x,e]$ is a non-zero element of $f_i \sG f_i$ with trace zero.  
The
 assumption implies that $[x,e] \in \sG \subseteq \sG + \sK$, as
required for $\sG + \sK$ to be a Lie ideal.

Now suppose that $\sG$ is a Lie addend for $\sK$.  Let $\sL$ 
denote
$\sG + \sK$.  While $f_i \sL f_i = f_i \sG f_i$ need not be a Lie
ideal in $\sL$ (indeed, need not even be a subset of $\sL$), it 
is
easy to see that $f_i \sL f_i$ is a Lie ideal in the full matrix
algebra $f_i \sB(\sH) f_i$.  Since there are only four Lie 
ideals in a
full matrix algebra, this shows that $f_i \sG f_i$ is one of the 
four
subspaces cited in Proposition \ref{Lieadd_struct}.

Suppose that $f_i g f_i$ is an element of $f_i \sG f_i$ with 
trace
zero.  If $f_i g f_i =0$, then it is certainly an element of 
$\sG$.  
If  $f_i g f_i \neq 0$, then there is an element $h$ in
$f_i \sB(\sH) f_i$ which does not commute with $f_i g f_i$.  
But $h \in \sA$ and $[g,h] = [f_i g f_i, h]$ is then a non-zero
element of $\sL$.  Let $x$ denote this element.  Observe that 
$x$ is a
non-zero element of $f_i \sG f_i$ which has trace zero.  But $x$ 
is
also in $\sL$; therefore, the smallest Lie ideal containing $x$ 
in the
algebra $f_i \sB(\sH) f_i$ is also contained in $\sL$.  This is 
the
Lie algebra of trace zero matrices in $f_i \sB(\sH) f_i$.  Since 
this
is also a subset of $\sE$, it is a subset of $\sG$.  Thus, all
elements of $f_i \sG f_i$ with trace zero are elements of $\sG$.

All that remains is to prove the constraint conditions when
$(i,j) \in \tilde S \setminus \tilde K$.  The following Lemma
decreases the need for cumbersome notation.

\begin{lemma} \label{two_summand}
Let $\sH = \sH_1 \oplus \sH_2$ be a direct sum of Hilbert spaces 
and
let $f_1$ and $f_2$ be the orthogonal projections on $\sH_1$ and
$\sH_2$. Let $x \in f_1 \sB(\sH)f_1 + f_2 \sB(\sH)f_2$.  Assume 
that
$[x,d]=0$ for all $d \in \sB(\sH)$ for which $d = f_1df_2$.  Then
there is a scalar $\lambda$ such that $x = \lambda f_1 + \lambda 
f_2$.
\end{lemma}

\begin{remark} 
Of course, the conclusion simply says that $x$ is a scalar 
operator.
But in the application of the Lemma, there will be
 additional summands present, so the
$\lambda f_1 + \lambda f_2$ format is more suitable.
\end{remark}

\begin{proof}
Let $\alpha_1, \alpha_2, \dots$ be an orthonormal basis for 
$\sH_1$
and $\beta_1, \beta_2, \dots$ an orthonormal basis for $\sH_2$.  
Let 
$d = \alpha_i \beta_j^*$ be the rank one partial isometry in
$\sB(\sH)$ with initial space $\bbC \beta_j$ and final space
$\bbC \alpha_i$.  Since $d = f_1 d f_2$, 
\[
f_1 x f_1 d - d f_2 x f_2 = [f_1xf_1+f_2xf_2,d]=[x,d]=0.
\]
For any $p$, 
$\langle f_1xf_1d\beta_j,\alpha_p \rangle = 
\langle f_1xf_1\alpha_i, \alpha_p \rangle =
\langle x\alpha_i, \alpha_p \rangle$.  If $p \ne i$ then
$\langle df_2xf_2\beta_j, \alpha_p \rangle = 0$.  (The range of
$d$ is $\bbC \alpha_i$, which is orthogonal to $\alpha_p$.)  It
follows that
\[
\langle x\alpha_i, \alpha_p \rangle =
\langle (f_1xf_1d-df_2xf_2)\beta_j, \alpha_p \rangle = 0.
\]
This is valid for all pairs of indices $p$ and $i$ with $p \ne 
i$.

Similarly, for any $q$,
$\langle df_2xf_2\beta_q, \alpha_i \rangle =
\langle f_2xf_2 \beta_q, \beta_j\rangle =
\langle x\beta_q, \beta_j \rangle$.  And if $q \ne j$, then
$f_1xf_1d\beta_q = 0$.  Hence
\[
\langle x\beta_q, \beta_j \rangle =
- \langle (f_1xf_1d-df_2xf_2)\beta_q, \alpha_i \rangle = 0.
\]
This too is valid for all pairs $q$ and $j$ with $q \ne j$.

We also have
\begin{align*}
\langle f_1xf_1d\beta_j, \alpha_i \rangle &=
\langle x\alpha_i, \alpha_i \rangle , \text{ and}\\
\langle df_2xf_2\beta_j, \alpha_i \rangle &=
\langle x\beta_j, \beta_j \rangle.
\end{align*}
Therefore,
\[
\langle x\alpha_i, \alpha_i \rangle - \langle x\beta_j, \beta_j
\rangle = \langle (f_1xf_1d-df_2xf_2) \beta_j, \alpha_i \rangle 
= 0
\]
and 
$\langle x\alpha_i,\alpha_i\rangle=\langle x 
\beta_j,\beta_j\rangle$.
This holds for any pair $i$ and $j$.  Letting $\lambda$ be this 
common
value, we now have
$x = \lambda f_1 + \lambda f_2$.
\end{proof}

It remains to show that if $f_i \sA f_j \ne \{0\}$ and
$f_i \sL f_j = \{0\}$, then for each $x \in \sG$ there is a 
scalar
$\lambda$ such that $f_ixf_i=\lambda f_i$ and 
$f_jxf_j = \lambda f_j$.  Recall that when $f_i \sA f_j \ne 
\{0\}$,
then
$f_i \sA f_j = f_i \sB(\sH) f_j$.  Let 
$d \in \sB(\sH)$ be such that $d = f_i d f_j$.  For $x \in \sG$,
$[x,d] = f_i[x,d]f_j \in f_i \sL f_j$, so $[x,d]=0$.  An 
application
of Lemma~\ref{two_summand} completes the proof of the Proposition
\ref{Lieadd_struct}.
\end{proof}

We now have a complete description of the structure of a Lie 
ideal in
a digraph algebra.
  Using this description, we give an alternate proof of
Theorem \ref{Th_sim} in the digraph algebra context:

\begin{proof}[Alternate proof of Theorem \ref{Th_sim} for digraph
  algebras]
We only need to show that if $\sL$ is a Lie ideal then it is 
invariant
under all similarity transforms.  Write $\sL$ in the form
$\sG + \sK$ as above.  Since $\sK$ is an associative ideal, it is
invariant under any similarity transform.  This reduces the 
proof to
showing that if $x \in \sG$ and if $t \in \sA$ is invertible, 
then
$t^{-1}xt \in \sG + \sK$.

If $t$ is an invertible element of $\sA$, then its diagonal part
$d=\sum_i f_itf_i$ is an invertible element of $\sE$.  (The 
inverse is
$\sum_i f_i t^{-1} f_i$.)  Now $d^{-1}t$ is invertible in $\sA$ 
and
can be written in the form
$d^{-1}t=1+n$, where $1$ is the identity element and $n$ is 
strictly
block upper triangular and therefore nilpotent.  So we can split 
the
argument into two cases.

\emph{Case 1}:  Assume that $x \in \sG$ and that $d$ is an 
invertible
element of $\sE$.  For each $i$, let $x_i = f_ixf_i$ and
$d_i = f_idf_i$.  Then
$d^{-1}xd = \sum_i d_i^{-1}x_id_i$. (Here, $d_i^{-1}$ is to be
interpreted as the inverse of $d_i$ in $f_i \sB(\sH) f_i$.)
 If $i$ is such that $f_i\sG f_i$ is either $\{0\}$ or the trace 
zero
 elements of $f_i \sB(\sH) f_i$, then $d_i^{-1}xd_i$ is a trace 
zero
 element of $f_i \sG f_i$ and so is in $\sG$.
 If $i$ is such that $f_i \sG f_i$ are scalar
elements, then $d_i^{-1}x_id_i = x_i$; in particular,  all the 
scalar
constraints are preserved by conjugation by $d$. 
This leaves $\sum d_i^{-1}xd_i$ where the sum is taken over 
those $i$
for which $f_i \sG f_i$ is a full matrix algebra.  There may be
constraints involving these indices, but they are all on
$\tr x_i$ only.  Since similarity preserves traces, the 
$d_i^{-1} x d_i$ satisfy the same constraints, so
$\sum d_i^{-1}xd_i \in \sG$.
 Thus
$d^{-1}xd$ is an element of $\sG \subseteq \sG + 
\sK$.\footnote{We
 do not need the details, but
  the constraints which determine $\sG$ as a subspace of $\sF$ 
take
  the form of a set of linear equations in the indeterminates 
$T_i$,
  where the $T_i$ correspond to those $f_i$ for which $f_i \sG 
f_i$ is
  either the full matrix algebra or the scalars.  The elements of
  $\sG$ are those $x \in \sF$ such that the numbers
 $\tr (f_i^{-1} x f_i)$ satisfies all the equations.} 

\emph{Case 2}: Assume that $x \in \sG$ and that $n$ is a 
nilpotent
element of $\sA$.  Let $k$ be the order of nilpotence for $n$.
Observe that
\begin{gather*}
(1+n)^{-1}x(1+n) = (1-n+n^2-n^3+\dots+(-1)^kn^k)x(1+n) \\
  =x+(-nx+xn)+(-nxn+n^2x)+(n^2xn-n^3x) \\
 +(-n^3xn+n^4x)
+\dots+((-1)^{k-1}n^{k-1}xn+(-1)^kn^kx) + (-1)^kn^kxn \\
=x-[n,x]+n[n,x]-n^2[n,x]+\dots+(-1)^{k+1}n^{k}[n,x].
\end{gather*}
The last equality uses the fact that $n^{k+1}=0$.
Since $x \in \sG \subseteq \sL$, $[n,x] \in \sL$.  But $n$ is 
strictly
block upper triangular in $\sA$; it follows that $[n,x]$ is also
strictly block upper triangular.  Hence, $[n,x] \in \sK$.  Since 
$\sK$
is an associative ideal, $n^j[n,x] \in \sK$ for all $j$.  This 
shows
that $(1+n)^{-1}x(1+n) \in \sG + \sK$, completing the proof of 
the
theorem. 
\end{proof}

\subsection{Triangular algebras (finite dimensional)}

Finite dimensional triangular algebras form a subclass of the 
digraph
algebras; consequently the work above on
digraph algebras  gives a description of
the Lie ideals in a finite dimensional triangular algebra.
  Since our goal is to extend these results to triangular
subalgebras of AF \cstar algebras, we pause to describe
explicitly the specialization of the digraph algebra results
 to the triangular algebra context.
There is more than one way to phrase this description; the
one used below was selected for its compatibility with the use of
groupoids in studying subalgebras of AF \cstar algebras.

Any finite dimensional triangular operator algebra is 
isomorphic to a subalgebra of the upper triangular
matrices $\sT_n$ which contains the diagonal $\sD_n$, for some
positive integer $n$.  Accordingly,  $\sB = M_n(\bbC)$ is the
context for the following discussion.

With $n$ fixed, we let $\sD$ be the algebra of diagonal matrices;
$\sT$, the algebra of upper triangular matrices; and $\sS$, the 
algebra
of strictly upper triangular matrices
in $\sB$.  If $\sC \subseteq \sB$ is a
bimodule over $\sD$, let 
\[
\spec (\sC) =\{(i,j) \mid c_{ij} \ne 0, \text{ for some }
c \in \sC\}
\]
It is easy to check that 
$\sC = \{x \in \sB \mid x_{ij}=0 \text{ whenever } (i,j) \notin
\spec(\sC)\}$.  

Let $\sA$ be a triangular subalgebra of $\sT$ (so  
$\sA \cap \sA^* = \sD$).
 Suppose that $\sK$ is an associative
ideal in $\sA$ and that 
$\sK \cap \sD = \{0\}$ (i.e., $\sK \subseteq \sA \cap \sS$).  
Then
 $\sK$ is the smallest Lie ideal in $\sA$ with the property
that $\sL \cap \sS = \sK$. Let
 $A = \spec(\sA)$, $D = \spec(\sD)$ and
$K = \spec(\sK)$.

Let
$\sE = \{d \in \sD \mid d_{ii} = d_{jj} \text{ whenever }
(i,j) \in A \setminus K\}$.  Then
$\sE+\sK$ is a Lie ideal and  is the largest Lie ideal such
that $\sL \cap \sS = \sK$.  Furthermore, if $\sF$ is a subspace 
of
$\sE$, then $\sF + \sK$ is a Lie ideal whose off-diagonal part is
$\sK$.  Letting
$\Lie(\sK) =\{\sF+\sK \mid \sF \subseteq \sE\}$, $\Lie(\sK)$ is 
the
family of all Lie ideals whose off-diagonal part is $\sK$.  
Finally,
 if $\sL$ is any Lie ideal in $\sA$, and if
$\sK = \sL \cap \sS$, then $\sK$ is an associative ideal and 
$\sL \in \Lie (\sK)$.

A direct proof of this description of the Lie ideals in a finite
dimensional triangular algebra is somewhat simpler than the
argument for the digraph algebra case.  
The primary simplification arises from the fact that all blocks 
are
$1 \times 1$ and therefore have no non-zero elements with trace 
zero.
For the full upper triangular
matrix algebra case, the description is covered by the 
literature on
Lie ideals in nest algebras. 

\section{Triangular subalgebras of AF \cstar algebras.} 
\label{AF}

In this section, $\sB$ will denote an AF \cstar algebra with 
canonical
diagonal $\sD$ and $\sA$ will be a 
 triangular subalgebra of $\sB$ with
diagonal $\sD$ (i.e., $\sA \cap \sA^* = \sD$).  This implies 
that there
is a sequence $\sB_n$ of finite dimensional \cstar algebras, 
each with
a maximal abelian self-adjoint subalgebra $\sD_n$ such that 
$\ds \sB = \lim_{\longrightarrow}\sB_n$ and
$\ds \sD = \lim_{\longrightarrow}\sD_n$.  We can, and do, view 
the
$\sB_n$ as a chain of subalgebras of $\sB$ in the usual way.
Since $\sA$ is a bimodule over $\sD$, $\sA$ is inductive.  This 
means
that $\sA_n \underset{\text{def}}{=} \sA \cap \sB_n$ is a 
triangular
subalgebra of $\sB_n$ with diagonal $\sD_n$ and 
$\ds \sA = \lim_{\longrightarrow}\sA_n$.

In addition to the presentation for $\sA$ described above, we 
shall
make use of ``coordinitization'' for $\sA$.  Coordinitization, or
groupoids, for \cstar algebras is treated in detail in the books 
of 
Renault \cite{MR82h:46075} and Paterson \cite{MR2001a:22003}.  
For a
good introduction to the use of groupoids in non-self-adjoint
algebras, see Muhly and Solel \cite{MR90m:46098}.  For the 
convenience
of the reader, we provide a brief sketch of the most relevant 
aspects
of coordinitization.

 Since $\sB$ is an AF
\cstar algebra, it is a groupoid \cstar algebra.  Let $G$ be the
groupoid.  The groupoid can be realized as a topological 
equivalence
relation on a compact Hausdorff, completely disconnected set 
$X$; $X$
will be such that $\sD \cong C(X)$. It is possible to pick a 
system of
matrix units $\{e^n_{ij}\}$ for $\sA$ such that, for each $n$,
 $\{e^n_{ij}\}$ are the matrix units which generate $\sA_n$ and 
each
 matrix unit in $\sA_n$ can be written as a sum of matrix units 
in
$\sA_{n+1}$.  The matrix units are all normalizing partial 
isometries
for $\sD$.  (A partial isometry $v$ is normalizing if
$v^*\sD v \subseteq \sD$ and $v\sD v^* \subseteq \sD$.)  The 
action of a
normalizing partial isometry on $\sD$ induces a partial 
homeomorphism
on $X$ and the equivalence relation $G$
is exactly the union of the graphs
of all the partial homeomorphisms induced by normalizing partial
isometries. 
The multiplication on $G$ is defined for those pairs of elements
$(x,y)$ and $(w,z)$ for which $w=y$; the product for a 
composable pair
is given by $(x,y)(y,z)=(x,z)$. Inversion is given by
$(x,y)^{-1} = (y,x)$.
 The topology on $G$ is the one obtained by declaring that
each such graph is an open subset of $G$.  It turns out that 
these
sets are all also closed, in fact, compact.  This description 
makes it
clear that the groupoid is independent of the presentation, but 
a very
handy fact is that $G$ is the union of the graphs of the matrix 
units
in the presentation.

 The elements of $\sB$ can be
identified with elements of $C_0(G)$ (but not all elements of 
$C_0(G)$ correspond to elements of $\sB$).  We won't need all the
details of this, but we will need the formula for multiplication:
if $f$ and $g$ are elements of $\sB$ viewed as functions in 
$C_0(G)$, then
\[
f\cdot g (x,y) = \sum_u f(x,u) g(u,y)
\]
where $u$ varies over the equivalence class of 
$x$ (which is the same as the equivalence class of $y$).
Note, in particular, that if $g \in \sD$ then the support of
$g$ is in $D$ and $f \cdot g (x,y) = f(x,y)g(y,y)$ and
$g\cdot f (x,y) = g(x,x)f(x,y)$.

If $\sC \subseteq \sB$ is any bimodule over $\sD$, then 
$C = \{(x,y)\in G \mid f(x,y) \ne 0 \text{ for some } f \in 
\sC\}$ has
the property that
$\sC = \{f \in \sB \mid \supp(f) \subseteq C\}$.  (This is the
spectral theorem for bimodules  \cite{MR90m:46098}; 
we shall refer to $C$ as the spectrum
of $\sC$ and  write $C = \spec(\sC)$.)  

With this terminology, 
$\spec(\sD) = \{(x,x) \mid x \in X\} $.  It is customary to 
identify 
$D=\spec(\sD)$ with $X$ (which is the spectrum of $\sD$ in the 
usual
sense for abelian \cstar algebras) by writing $x$ in place of 
$(x,x)$,
but we won't do so in this treatment.
If $A = \spec(\sA)$, then $A$ is a subrelation of $G$ whose
intersection with its reversal (the spectrum of $\sA^*$) is 
exactly
$D$.  Let $S = A \setminus D$.  If $f \in \sA$, then we can write
$f = f|_D + f|_S$; this gives a decomposition of $f$ into a 
diagonal
part and an off-diagonal part.

There is another way to effect the same decomposition, via a
contractive conditional expectation onto the diagonal.  To 
define this
conditional expectation we use the presentation 
$\ds \lim_{\longrightarrow}\sB_n$ and the matrix unit system 
$\{e^n_{ij}\}$.  For each $n$, define $\pi_n$ on $\sB$ by
$\pi_n(f) = \sum_i e^n_{ii}fe^n_{ii}$.  The sequence of maps
$\pi_n$ converges pointwise to a contractive conditional 
expectation
of $\sB$ onto $\sD$.  If $f \in \sA$, it can be shown that
$\pi(f) = f|_D$ and $f-\pi(f) = f|_S$.  In particular, the 
conditional
expectation is independent of the choice of presentation or the 
choice
of matrix units.

If $f$ and $g$ are elements of $\sA$, then the index of
summation in the formula for the
product $f\cdot g (x,y) = \sum_u f(x,u) g(u,y)$, 
runs over all elements  $u$ such that
$(x,u) \in A$ and $(u,y) \in A$. (This is always a countable
set.)  In particular, if $x=y$ then the
only possible value for $u$ is $x$, so
$f \cdot g (x,x) = f(x,x)g(x,x)$.  From this it immediately 
follows
that all commutators vanish on $D$.  

Let $\sK$ be an associative
ideal in $\sA$ such that $\sK \cap \sD = \{0\}$.  
Let $K = \spec(\sK)$. Thus $K \cap D = \emptyset$ and
$K \subseteq S$.  Also, $\pi(f) = 0$ for all $f \in \sK$.  
Trivially,
$\sK$ is a Lie ideal in $\sA$.  Now let
$\sE = \{f \in \sD \mid f(x,x)=f(y,y) \text{ whenever } (x,y) 
\in S
\setminus K\}$.  Let $\sF$ be any subspace of $\sE$ and let
$\sL = \sF + \sK$.  With this notation:

\begin{proposition} \label{maxideal}
$\sL$ is a Lie ideal in $\sA$.
\end{proposition}

\begin{proof}
  It suffices to prove that
$[f,e] \in \sL$ for any $f \in \sF$ and any matrix unit $e$ in a
matrix unit system for $\sA$.  If $e$ is a diagonal matrix unit,
then $[f,e]=0 \in \sL$. So assume that $e$ is off-diagonal.  
Viewed as
a function in $C_0(G)$, $e$ is the characteristic function of a
$G$-set which is contained in $S$.  Since the support of $[f,e]$ 
is
disjoint from $D$, it is contained in $S$.  If we can show that 
it is
contained in $K$, then $[f,e] \in \sK \subseteq \sL$.  So let
$(x,y) \in S \setminus K$.  Then
\begin{align*}
[f,e](x,y) &=
f\cdot e (x,y) - e\cdot f (x,y) = f(x,x)e(x,y)-e(x,y)f(y,y) \\
&=(f(x,x)-f(y,y))e(x,y) = 0.
\end{align*}
 
This shows that $[f,e]$ is supported in $K$.
\end{proof}

The next Theorem shows that all closed Lie ideals have this form.

\begin{theorem} \label{Lieform}
If $\sL$ is a Lie ideal in $\sA$ then $\sL = \sF + \sK$,
 where $\sK$ is a diagonal disjoint associative ideal in
$\sA$ and $\sF$ is a subspace of 
\[
\sE_{\sK} =
\{f \in \sD \mid f(x,x)=f(y,y) \text{ whenever } (x,y) \in S
\setminus K\}.
\]
\end{theorem}

\begin{proof}
Assume that $\sL$ is a Lie ideal in $\sA$.  We may define $\sK$ 
in
any of several equivalent ways:
\begin{align*}
\sK &= \{f \in \sL \mid \pi(f) = 0 \} \\
&= \{f-\pi(f) \mid f \in \sL\} \\
&= \{f \in \sL \mid \supp(f) \subseteq S\}.
\end{align*}
We must show that $\sK$ is an associative
ideal in $\sA$; the first step in that
direction is to show that $\sK$ is a bimodule over $\sD$.  (The
inductivity of $\sK$ will be helpful in showing that $\sK$ is an
associative ideal.)
The following lemma is useful in showing that $\sK$ is a bimodule
over $\sD$.

\begin{lemma} \label{row}
Let $\sL$ be a Lie ideal in $\sA$ and let $f \in \sL$.  Let
$d_1, \dots, d_q$ be the minimal diagonal projections in
$\sD_m$.  Then, for each $i$, 
 $d_if-d_ifd_i$ is an element of $\sL$.
\end{lemma}

\begin{proof}
Fix $i$. If $j \neq i$ then
$[[d_i,f],d_j]=[d_i f - f d_i, d_j] = d_i f d_j + d_j f d_i 
\in \sL$.  Hence,
\[
\sum_{j \ne i}  d_ifd_j+ \sum_{j \ne i}d_jfd_i \in \sL.  
\]
We also have
\begin{align*}
[d_i,f] &= d_if-fd_i \\
&= \sum_{j=1}^n d_ifd_j - \sum_{j=1}^n d_jfd_i \\
&= \sum_{j \ne i} d_ifd_j - \sum_{j \ne i} d_j fd_i \in \sL.
\end{align*} 
 Take an average to see that
$\ds \sum_{j \ne i} d_i f d_j = d_i f (1-d_i) = d_if-d_ifd_i \in 
\sL$. 
(Or, more succinctly,
$\ds d_if-d_ifd_i = \frac 12 ([d_i,f]+\sum_{j\ne i}[d_i,[f,d_j]] 
\in \sL$.)
\end{proof}

Let $f \in \sK$.  Fix $n$ and let $e_1, \dots, e_p$ be the 
minimal
diagonal projections in $\sD_n$.  As above
\[
f-\pi_n(f) = \sum_{i \ne j}e_ife_j.
\]
and $e_ife_j + e_jfe_i \in \sL$ for each pair
$i,j$ with $i \ne j$.  By the way, even when 
$e_j \sA_n e_i = \{0\}$, it is possible that both
$e_ife_j$ and $e_jfe_i$ are non-zero.

In order to show that $\sK$ is a left
$\sD$ module,  it suffices to show that
 $d(e_ife_j +e_jfe_i)\in \sK$ for all 
$d \in \sD$ and all pairs ${i,j}$ with $i \ne j$.
 It then follows that
$df-d\pi_n(f) \in \sK$ and, since $\sK$ is closed,
 $df = \lim_n (df -d\pi_n(f)) \in \sK$.  
 To do this for all $d$, it is enough to show that
$d(e_ife_j+e_jfe_i) \in \sK$ when 
$d$ is a minimal  diagonal matrix unit in $\sD_m$ and
$m \geq n$.

Fix $m \geq n$ and let $d$ be a minimal diagonal matrix unit in 
$\sD_m$.  
Since $m \geq n$, $d$ is a subprojection of one of the 
$e_k$.
In particular, since $i \ne j$, either $de_i=0$ or $de_j=0$; in 
either
case $d(e_ife_j+e_jfe_i)d=0$.  But then Lemma \ref{row} yields
\[
d(e_ife_j+e_jfe_i) = d(e_ife_j+e_jfe_i)-d(e_ife_j+e_jfe_i)d 
\in \sL.
\] 
 But
$\pi(d(e_ife_j+e_jfe_i)) = d\pi(e_ife_j+e_jfe_i)=0$, so
$d(e_ife_j+e_jfe_i) \in \sK$, as desired.

This shows that $\sK$ is a left $\sD$-module.  A similar argument
shows that $\sK$ is a right $\sD$-module.  Or, alternatively, 
$\sK^*$
is  the off-diagonal part of a closed Lie ideal in $\sA^*$
 and so is a left
$\sD$-module, which immediately implies that $\sK$ is a right 
$\sD$-module.  One consequence of this fact is that we now know 
that
$\sK$ is the closed linear span of all the off-diagonal matrix 
units
in $\sL$.  (Note: the diagonal part of $\sL$ is not inductive, in
general.) 

Next, we show that $\sK$ is a left ideal.  Let $f \in \sK$.  It
suffices to show that $ef \in \sK$ for all matrix units $e$ in
$\sA$.  Assume that $e$ is an off-diagonal matrix unit in 
$\sA_n$.
(We do not need to consider diagonal matrix units, since we know 
that
$\sK$ is a $\sD$-module.)

Let $\epsilon >0$.  We claim that there is $g \in \sK$ such that
$\|g-ef\|<\epsilon$, i.e., $\dist (ef,\sK) < \epsilon$.  Since 
$\sK$
is closed and $\epsilon$ is arbitrary, the claim implies that
$ef \in \sK$ and $\sK$ is a left ideal.

Let $p$ be a positive integer such that $p \geq n$ and
$\|\pi_q(f)\| < \epsilon/3$ for all $q \geq p$.  (We can find 
such $p$
since $\pi_p(f) \rightarrow \pi(f)=0$ as 
$p \rightarrow \infty$.)  With $q \geq p$, let
$b \in \sK \cap \sB_q = \sK \cap \sA_q$ be
 such that $\|b-f\|<\epsilon/3$.  ($\sK$ is
inductive.) Let $c = b-\pi_q(b)$.  Since $b \in \sK$ and
$\sK$ is a bimodule over $\sD$, 
$\pi_q(b) \in \sK \subseteq \sL$.  Therefore 
 $c \in \sK \cap \sA_q \subseteq \sL \cap \sA_q$.  Now
\begin{align*}
\|c-f\| &= \|b-\pi_q(b)-f\| \\
&=\|b-f-\pi_q(f) + \pi_q(f-b)\| \\
&\leq \|b-f\| + \|\pi_q(f)\| + \|\pi_q(f-b)\| \\
&= \epsilon/3 +  \epsilon/3 + \epsilon/3 =\epsilon.
\end{align*}
Thus, $\|ec-ef\| \leq \|c-f\| \leq \epsilon$. So  we just need to
show that $ec \in \sK$.

Since $e$ is an off-diagonal matrix unit in $\sA_n$ and $q \geq 
n$, 
$e$ can be written as a sum of off-diagonal matrix units in 
$\sA_q$.  So we need to prove that $hc \in \sK$ when $h$ is an
off-diagonal matrix unit in $\sA_q$, and for this it is enough to
prove that $hc \in \sL$, since $\pi(hc)=0$.  Let
$\sS_q$ be the subalgebra of $\sA_q$ generated by all the 
off-diagonal
matrix units in $\sA_q$ and let
$\sL'=\sL \cap \sS_q$.  If $x \in \sA_q$ and 
$y \in \sL' \subseteq \sL$, then $[x,y]\in \sL$, since $\sL$ is 
a Lie
ideal. But $[x,y] \in \sS_q$ also, since it is a commutator of 
two
elements in $\sA_q$.  Thus $\sL'$ is a diagonal disjoint Lie 
ideal in
$\sA_q$.  But this implies that $\sL'$ is an associative ideal
in $\sA_q$ (Proposition \ref{dig_Lieid_struct}).
Since $c \in \sL \cap \sA_q$ and $\pi_q(c)=0$, $c \in \sL'$.
Therefore,
$hc \in \sL' \subseteq \sL$.
This completes the argument that $\sK$ is a left ideal.  It also 
shows
that $\sK^*$ is a left ideal in $\sA^*$,
 whence $\sK$ is a right ideal in $\sA$.

We have now shown that  if $\sL$ is a Lie ideal in $\sA$ then the
diagonal disjoint part $\sK$ is an 
associative ideal in $\sA$.  Let 
$\sF = \sL \cap \sD$.  From the way that $\sK$ is defined, it is 
clear
that
$\sL = \sF + \sK$.  To complete the description of $\sL$, we 
need to
show that if $f \in \sF$ then $f(x,x)=f(y,y)$ whenever
$(x,y) \in S\setminus K$.  Let $(x,y) \in S \setminus K$.  Let 
$e$ be
a matrix unit in $\sA$ such that $(x,y) \in \supp (e)$.  Since 
$f \in \sL$, the commutator $[f,e]$ is also in $\sL$. 
Since commutators vanish on the diagonal,
$[f,e]\in \sK$.  Hence
\begin{gather*}
0 = [f,e](x,y) = f \cdot e (x,y) - e \cdot f (x,y) \\
= f(x,x)e(x,y)-e(x,y)f(y,y) = f(x,x)-f(y,y).
\end{gather*}  
Thus,
$f(x,x)=f(y,y)$ whenever $(x,y) \in S \setminus K$.
\end{proof}

With $\sL$ and $\sK$ as above, if we let 
\[
\sE_{\sK} = \{f \in \sD \mid f(x,x)=f(y,y)
 \text{ whenever } (x,y) \in S
\setminus K \},
\]
 then we have shown that $\sK$ is the smallest Lie
ideal and $\sE_{\sK} + \sK$ is the largest Lie ideal 
 which has $\sK$ as its
off-diagonal part.
Once again, the structure of Lie ideals in triangular AF algebras
permits an alternative proof of Theorem \ref{Th_sim} for this
context. 

\begin{proof}[Alternate proof of Theorem \ref{Th_sim} for 
triangular
subalgebras of AF \cstar algebras]
As usual, we need only prove that closed Lie ideals are invariant
under similarities.
Let $\sL$ be a closed Lie ideal and write $\sL = \sF + \sK$ as 
above.
Let $t$ be an invertible element of $\sA$.  Since the invertible
elements of an operator algebra form an open set, there is a 
sequence
$t_p \in \sA_p$ of invertible elements in $\sA_p$ such that
$t_p \rightarrow t$ and $t_p^{-1} \rightarrow t^{-1}$.  If we
can show that $t_p^{-1} \sL t_p \subseteq \sL$ for all $p$, then
$t^{-1} \sL t \subseteq \sL$.  

So we have reduced the proof to the case where $t \in \sA_p$ for 
some
$p$. Since $\sA_p$ is isomorphic to a triangular matrix algebra, 
we
can write $t=d(1+n)$, where $d$ is a diagonal invertible element 
of
$\sA_p$ and $n$ is a nilpotent element of $\sA_p$.  Let $k$ be 
the
order of nilpotence of $n$.  

Since $\sK$ is an associative ideal, $t^{-1}\sK t \subseteq \sK$.
  So we need prove
that $t^{-1}\sF t \subseteq \sF + \sK = \sL$.  Since $d$ is 
diagonal
$d^{-1}\sF d = \sF$.  This leaves invariance under conjugation by
$1+n$.  For any $f \in \sF$,
\begin{gather*}
(1+n)^{-1}f(1+n) = (1-n+n^2-n^3+\dots+(-1)^kn^k)f(1+n) \\
  =f+(-nf+fn)+(-nfn+n^2f)+(n^2fn-n^3f) \\
 +(-n^3fn+n^4f)
+\dots+((-1)^{k-l}n^{k-1}fn+(-1)^kn^kf) +(-1)^kn^kfn \\
=f-[n,f]+n[n,f]-n^2[n,f]+\dots+(-1)^{k+1}n^{k}[n,f].
\end{gather*}
Now $[n,f] \in \sL$ and has off-diagonal support, 
so $[n,f] \in \sK$.
Since $\sK$ is an associative ideal, 
$(1+n)^{-1}f(1+n) \in \sF + \sK = \sL$.  Thus
$t^{-1}\sL t \subseteq \sL$; $\sL$ is invariant under inner 
automorphisms.
\end{proof}

%\bibliographystyle{amsplain}
%\bibliography{biblio}

\providecommand{\bysame}{\leavevmode\hbox 
to3em{\hrulefill}\thinspace}

\end{document}